\documentclass[11pt]{amsart}
\usepackage{latexsym,amssymb,amsmath,hyperref}
\usepackage{amsthm, dsfont}
\usepackage[all]{xy}
\usepackage[short,nodayofweek]{datetime}
\usepackage[active]{srcltx}

\leftmargin  \leftmargini
\def\ZZ{{\mathbb Z}}

\def\RR{{\mathbb R}}
\def\CC{{\mathbb C}}

\def\H{{\mathcal H}}

\DeclareMathOperator{\GL}{GL}

\DeclareMathOperator{\trace}{tr}

\DeclareMathOperator{\Sp}{Sp}

\theoremstyle{plain}
\newtheorem{thm}{Theorem}[section]
\newtheorem{cor}[thm]{Corollary}
\newtheorem{lem}[thm]{Lemma}

\newtheorem{eg}[thm]{Example}

\theoremstyle{definition}

\theoremstyle{remark}

\newtheoremstyle{Acknowledgements}
  {}
    {}
     {}
     {}
    {\bfseries}
    {}
     {.5em}
     {\thmname{#1}\thmnumber{ }\thmnote{ (#3)}}
\theoremstyle{Acknowledgements}

\date{\today, \currenttime} 

\begin{document}

\title[Casimir operators for symplectic groups]{Casimir operators for symplectic groups}

\author{Kathrin Maurischat}
\address{\rm {\bf Kathrin Maurischat}, Mathematics Center  Heidelberg
  (MATCH), Heidelberg University, Im Neuenheimer Feld 288, 69120 Heidelberg, Germany }
\curraddr{}
\email{\sf maurischat@mathi.uni-heidelberg.de}
\thanks{}




\begin{abstract}
We give a full set of Casimir operators for the symplectic group of arbitrary genus in terms of a basis  chosen such that the action on representations of known $K$-type becomes transparent. We give examples for the latter.
\end{abstract}

\maketitle




\section*{Introduction}
The original intension of this  work was to understand the action of Casimir operators on some kinds of automorphic forms for the symplectic group $\operatorname{Sp}_m(\RR)$.

In \cite{weissauersLN}, this problem is done for the standard first Casimir operator $C_1$. 
Using the Cartan decomposition of the symplectic lie algebra $\mathfrak g=\mathfrak k\oplus\mathfrak p$, one  decomposes $C_1=\trace(E_+E_-)+k_1$, where $\trace(E_+E_-)$ is a differential operator on the Siegel halfplane depending only on $\mathfrak p^\CC$ and $k_1$ is some constant coming from $\mathfrak k^\CC$ depending on the $K$-type only.
Surprisingly,  an analog for higher Casimir operators up to the genus $m$ does not exist in literature.
Usually Casimir operators are realized with respect to a Cartan subalgebra  which evidently is not of any help here. 

We use a basis of $\mathfrak g^\CC=\mathfrak k^\CC\oplus\mathfrak p_+\oplus\mathfrak p_-$ which has  pleasing properties: Lie multiplication as well as matrix multiplication is  simple and the dual basis (with respect to half the trace) is essentially deduced by rearranging.  The basis differs from that used in \cite{weissauersLN} in the $\mathfrak k^\CC$-part. For this basis we evaluate the common formula (\cite{knapp2}, IV.~7)
\begin{equation}\label{Casimir_allgemein}
D_r= \sum_{i_1,\dots,i_r}\trace(X_{i_1}\cdots X_{i_r})X_{i_1}^\ast\cdots X_{i_r}^\ast
\end{equation}
for Casimir elements to get a set $\{D_{2},\dots,D_{2m}\}$ of $m$ Casimir elements which indeed generates the center $\mathfrak Z(\mathfrak g^\CC)$ of the universal envelopping algebra.
As examples, we give precise formulae for $D_2, D_4$.

We apply the result to determine the action of $\mathfrak Z(\mathfrak g^\CC)$ on a representation of $K$-type $(\lambda,\dots,\lambda)$ to be given by that of $\trace(E_+E_-),\dots,\trace((E_+E_-)^m)$. For  automorphic forms,  the latter are differential operators on the Siegel halfplane.

\section{Notation}\label{section_notation}

Let $G=\operatorname{Sp}_m(\RR)$ be the real symplectic group of genus $m$ and let $\mathfrak g=\mathfrak{sp}_m(\RR)$ be its Lie algebra.
We consider the matrix realization of its complexification $\mathfrak g^\CC\subset M_{2m,2m}(\mathbb C)$ consisting of those $g$ satisfying
\begin{equation*}
 g'\begin{pmatrix}0&-\mathbf 1_m\\\mathbf 1_m&0\end{pmatrix}+\begin{pmatrix}0&-\mathbf 1_m\\\mathbf 1_m&0\end{pmatrix}g=0.
\end{equation*}
The Cartan decomposition for $\mathfrak g$ implies that
$\mathfrak g^\CC=\mathfrak k^\CC\oplus\mathfrak p_+\oplus\mathfrak p_-$,
where $\mathfrak k^\CC$ is given by those matrices satisfying
\begin{equation*}
 \begin{pmatrix}A&-S\\S&A\end{pmatrix}, \quad A'=-A, \quad S'=S,
\end{equation*}
and
\begin{equation*}
 \mathfrak p_\pm=\left\{\begin{pmatrix}X&\pm iX\\\pm iX&-X\end{pmatrix},\quad X'=X\right\}.
\end{equation*}
Let $e_{kl}\in M_{m,m}(\mathbb C)$ be the elementary matrix having entries $(e_{kl})_{ij}=\delta_{ik}\delta_{jl}$ and let $X^{(kl)}=\frac{1}{2}(e_{kl}+e_{lk})$. Further, let $A^{(kl)}=e_{kl}-e_{lk}$ and $S^{(kl)}=2X^{(kl)}$.

 The elements $E_{\pm kl}=E_{\pm lk}$ of $\mathfrak p_\pm$ are defined to be those corresponding to $X=X^{(kl)}$, $1\leq k,l\leq m$. 
Then $E_{\pm kl}$, $1\leq k\leq l\leq m$, form a basis of $\mathfrak p_\pm$.
For abbreviation, let $E_\pm$ be the matrix having entries $E_{\pm kl}$. 

Define the following elements of $\mathfrak k^\CC$  for $1\leq k,l\leq m$: Let $a_{kl}k^\CC$ be given by
$S=0$ and $A=A^{(kl)}$ and let $s_{kl}k^\CC$ be given by $A=0$ and $S=S^{(kl)}$.
Then a basis of $\mathfrak k^\CC$ is given by $B_{kl}:=\frac{1}{2}(a_{kl}+is_{kl})$, $1\leq k,l\leq m$.
Let $B=(B_{kl})_{kl}$ be the matrix with entries $B_{kl}$ and let $B^\ast$ be its transpose having entries $B_{kl}^\ast=B_{lk}$.

Lie multiplication in $\mathfrak g^\CC$ is easily checked to be given by
\begin{eqnarray*}
 &&[E_{+ij},E_{+kl}]=0,\quad\quad [E_{-ij},E_{-kl}]=0,\\
&&\left[ E_{+ij},E_{-kl} \right]=\delta_{ik}B_{jl}+\delta_{jl}B_{ik}+\delta_{il}B_{jk}+\delta_{jk}B_{il}\\
&&\left[ B_{ij},E_{+kl} \right]=\delta_{jk}E_{+il}+\delta_{jl}E_{+ik}, \\
&&[ B_{ij},E_{-kl} ]=-\delta_{ik}E_{-jl}-\delta_{il}E_{-jk}\\
&&\left[ B_{ij},B_{kl} \right] = \delta_{jk}B_{il}-\delta_{il}B_{kj}.
\end{eqnarray*}

We denote by $\mathcal B$ the nondegenerate bilinear form on $\mathfrak g^\CC$ defined by
\begin{equation}\label{def_bilinearform}
 \mathcal B(g,h)=\frac{1}{2}\trace(g\cdot h).
\end{equation}
With respect to $\mathcal B$ we get the following dual basis:
$E_{\pm kl}^\ast= \frac{1}{1+\delta_{kl}}E_{\mp kl}$ as well as
$B_{kl}^\ast=B_{lk}$ for all $k,l$.



\section{Casimir elements}\label{Section_casimire}
In the following, we study words in the matrices $E_+,E_-,B$ and $B^\ast$.
Let us define some conditions on these words:
\begin{itemize} 
\item[(i)] $E_+$ is followed by $E_-$ or $B^\ast$.
\item[(ii)] $E_-$ is followed by $E_+$ or $B$.
\item[(iii)] $B$ is followed by $E_+$ or $B$.
\item[(iv)] $B^\ast$ is followed by $E_-$ or $B^\ast$.
\item[(v)] $E_+$ occurs with the same multiplicity as $E_-$.
\end{itemize}
We start with a combinatorial lemma.
\begin{lem}\label{Lemma_woerter}
 Let $r>0$ be an integer. Then there are $2^{2r}$ possibilities to choose a word $w$ of length $2r$
in the matrices $E_+,E_-,B$ and $B^\ast$ such that the  conditions (i) to (v) are satisfied.
\end{lem}
\begin{proof}[Proof of Lemma~\ref{Lemma_woerter}]
We make a second claim slightly modifying Lemma~\ref{Lemma_woerter} and prove it along with the lemma itself by induction on $r$.

{\it Claim:} There are $2^{2r}$ possibilities to choose a word $w$ of length $2r$
in the matrices $E_+,E_-,B$ and $B^\ast$ such that the  conditions (i) to (iv)  and
\begin{itemize}
 \item [(v')] The multiplicity of $E_+$ is that of $E_-$ enlarged or reduced by one.
\end{itemize}
are satisfied.

For $r=1$, the four possible words of the lemma are $E_+E_-$, $E_-E_+$, $BB$ and $B^\ast B^\ast$, while the four possibilities of the claim are $E_+B^\ast$, $E_-B$, $BE_+$ and $B^\ast E_-$.
Now look at a word $w$  of length $2(r+1)$ satisfying (i)--(v). 
First, if $w$ ends with $E_+E_-$ or with $E_-E_+$, then the initial subword of length $2r$ satisfies (i)--(v). And for any of these initial subwords, the ending among $E_+E_-$ and $E_-E_+$ is unique. This gives $2^{2r}$ possibilities for $w$, using the lemma for $r$. Similarly, we get  $2^{2r}$ possibilities for a word where $E_\pm$ does not occur in the last two letters.
 If exactly one of the last two letters is $E_\pm$, then we again get $2\cdot 2^{2r}$ possibilities by   the claim for $r$. Together these are $2^{2(r+1)}$ possibilities. Similarly, we get the result for the claim, too.
\end{proof}
In the following, we  formally take the trace of a word $w$ in the operator valued matrices. For example,
\begin{equation*}
 \trace(E_+E_-)=\sum_{k,l}E_{+kl}E_{-lk}.
\end{equation*}

\begin{thm}\label{Satz_Casimir_allgemein}
Let $\mathfrak g=\mathfrak{sp}_m(\RR)$ be the Lie algebra of the symplectic group of genus $m$.

\begin{enumerate}
 \item [(a)]
 The $r$-th Casimir element is given by
\begin{equation*}
 D_{2r}=\sum_{w} (-1)^{L(w)}\trace(w),
\end{equation*}
where the sum is over all words $w$ of length $2r$ satisfying  conditions (i) to (v) above, and
 $L(w)$ is the number of times  $E_-B$  and $BE_+$ occur isolatedly   in $w$  counted cyclicly.
\item[(b)]
The center $\mathfrak Z(\mathfrak g^\CC)$ of the universal envelopping algebra of $\mathfrak g^\CC$ is generated by  the $m$ Casimir operators $D_2,\dots,D_{2m}$.
\end{enumerate}
\end{thm}
Here isolated means that $E_-B$ and $BE_+$ must not hit each other,  for example $L(E_-BE_+B^\ast)=1$ while $L(E_-BBE_+)=2$. And cyclic means that we have to take into account that the trace is cyclicly invariant, so e.g. $L(E_+E_-BB)=L(E_-BBE_+)=2$.
\begin{eg}\label{Prop_Casimir_2}
By Lemma~\ref{Lemma_woerter}, we have to sum over the traces of $2^{2r}$ words $w$. So the first two Casimirs are
\begin{equation*}
 D_2=\trace(E_+E_-)+\trace(E_-E_+)+\trace(BB)+\trace(B^\ast B^\ast),
\end{equation*}
\begin{eqnarray*}
 D_4 &=& \trace(E_+E_-E_+E_-)+\trace(E_-E_+E_-E_+)+\trace(BBBB)+\trace(B^\ast B^\ast B^\ast B^\ast)\\
&& +\sum_{\zeta\in Z_4}\bigl(\trace(\zeta(E_+E_-BB))+\trace(\zeta(E_-E_+B^\ast B^\ast))-\trace(\zeta(E_+B^\ast E_-B))\bigr),
\end{eqnarray*}
where $Z_4$ is the group of cyclic permutations of four elements.
\end{eg}

\begin{proof}[Proof of Theorem~\ref{Satz_Casimir_allgemein}]
(a) We define the following matrices
\begin{eqnarray*}
 K_1=\begin{pmatrix} \mathbf 1_m&i\mathbf 1_m\\-i\mathbf 1_m&\mathbf 1_m\end{pmatrix}, &&
K_2=\begin{pmatrix} \mathbf 1_m&-i\mathbf 1_m\\i\mathbf 1_m&\mathbf 1_m\end{pmatrix},\\
P_+=\begin{pmatrix} \mathbf 1_m&i\mathbf 1_m\\i\mathbf 1_m&-\mathbf 1_m\end{pmatrix},&&
P_-=\begin{pmatrix} \mathbf 1_m&-i\mathbf 1_m\\-i\mathbf 1_m&-\mathbf 1_m\end{pmatrix}.
\end{eqnarray*}
Notice that $K_j^2=K_j$, while $P_\pm^2=0$.
 In the following, if we abbriviate $e_{jk}K_{1/2}$, $e_{jk}P_{\pm}$, for a $(m\times m)$-elementary matrix $e_{jk}$ and a $(2m\times 2m)$-matrix $K_{1/2}, P_\pm$, we mean
\begin{equation*}
 e_{jk}K_1=\begin{pmatrix} e_{jk}&ie_{jk}\\-ie_{jk}&e_{jk}\end{pmatrix},\quad\textrm{ etc.}
\end{equation*}
Now we show that $D_{2r}$ in (\ref{Casimir_allgemein}) has the claimed shape. A single summand of $D_{2r}$ looks like
\begin{equation*}
 \trace(X_{j_1j_2}^{(1)}X_{k_1k_2}^{(2)}\cdots X_{j_{2r-1}j_{2r}}^{(2r-1)}X_{k_{2r-1}k_{2r}}^{(2r)})(X_{j_1j_2}^{(1)}\cdots X_{k_{2r-1}k_{2r}}^{(2r)})^\ast,
\end{equation*}
where $X_{j_{n}j_{n+1}}^{(n)}$ runs through the basis $E_{\pm jk}$, $1\leq j\leq k\leq m$, $B_{jk}$, $1\leq j,k\leq m$.
First we examine conditions for a pair $X_{j_1j_2}^{(n)}X_{k_1k_2}^{(n+1)}$ to occur in some summand. To get nice formulae, we sum over all pairs of the same kind.
First, let $X_{j_1j_2}^{(n)}=E_{+j_1j_2}$ and $X_{k_1k_2}^{(n+1)}=E_{-k_1k_2}$.
Computing the matrix product $E_{+j_1j_2}E_{-k_1k_2}$, taking duals and rearranging summation, we get
\begin{eqnarray}
 &&\hspace*{-10mm}\sum_{j_1\leq j_2;k_1\leq k_2}\frac{\trace(\cdots E_{+j_1j_2}E_{-k_1k_2}\cdots)}{(1+\delta_{j_1j_2})(1+\delta_{k_1k_2})}(\cdots)^\ast E_{-j_1j_2}E_{+k_1k_2}(\cdots)^\ast\label{gleichung_a}\\
&&=
\frac{1}{4}\sum_{j_1,j_2,k_1,k_2}\trace\bigl((\delta_{j_1k_1}e_{j_2k_2}+\delta_{j_1k_2}e_{j_2k_1}+\delta_{j_2k_1}e_{j_1k_2}+\delta_{j_2k_2}e_{j_1k_1})K_2\bigr)\nonumber\\
&&\hspace*{7cm}\cdot(\cdots)^\ast E_{-j_1j_2}E_{+k_1k_2}(\cdots)^\ast\nonumber\\
&&=
\sum_{j_1,j_2,k_1}\trace(\cdots e_{j_1k_1}K_2\cdots)(\cdots)^\ast E_{-j_1j_2}E_{+j_2k_1}(\cdots)^\ast.\nonumber
\end{eqnarray}
Similarly we get for the other choices of basis elements
\begin{eqnarray}
&&\hspace*{-10mm}\sum_{j_1\leq j_2;k_1\leq k_2}\trace(\cdots E_{-j_1j_2}E_{+k_1k_2}\cdots)(\cdots E_{-j_1j_2}E_{+k_1k_2}\cdots)^\ast\label{gleichung_b}\\
&&=
\sum_{j_1,j_2,k_1}\trace(\cdots e_{j_1k_1}K_1\cdots)(\cdots)^\ast E_{+j_1j_2}E_{-j_2k_1}(\cdots)^\ast,\nonumber
\end{eqnarray}
\begin{eqnarray}
&&\hspace*{-10mm}\sum_{j_1\leq j_2;k_1, k_2}\trace(\cdots E_{+j_1j_2}B_{k_1k_2}\cdots)(\cdots E_{+j_1j_2}B_{k_1k_2}\cdots)^\ast\label{gleichung_c}\\
&&=
-\sum_{j_1,j_2,k_1}\trace(\cdots e_{j_1k_1}P_+\cdots)(\cdots)^\ast E_{-j_1j_2}B_{j_2k_1}(\cdots)^\ast,\nonumber
\end{eqnarray}
\begin{eqnarray}
&&\hspace*{-10mm}\sum_{j_1\leq j_2;k_1, k_2}\trace(\cdots E_{-j_1j_2}B_{k_1k_2}\cdots)(\cdots E_{-j_1j_2}B_{k_1k_2}\cdots)^\ast\label{gleichung_d}\\
&&=
\sum_{j_1,j_2,k_1}\trace(\cdots e_{j_1k_1}P_-\cdots)(\cdots)^\ast E_{+j_1j_2}B^\ast_{j_2k_1}(\cdots)^\ast,\nonumber
\end{eqnarray}
\begin{eqnarray}
&&\hspace*{-10mm}\sum_{j_1, j_2;k_1\leq k_2}\trace(\cdots B_{j_1j_2}E_{-k_1k_2}\cdots)(\cdots B_{j_1j_2}E_{-k_1k_2}\cdots)^\ast\label{gleichung_e}\\
&&=
-\sum_{j_1,j_2,k_1}\trace(\cdots e_{j_1k_1}P_-\cdots)(\cdots)^\ast B_{j_1j_2}E_{+j_2k_1}(\cdots)^\ast,\nonumber
\end{eqnarray}
\begin{eqnarray}
&&\hspace*{-10mm}\sum_{j_1, j_2;k_1\leq k_2}\trace(\cdots B_{j_1j_2}E_{+k_1k_2}\cdots)(\cdots B_{j_1j_2}E_{+k_1k_2}\cdots)^\ast\label{gleichung_f}\\
&&=
\sum_{j_1,j_2,k_1}\trace(\cdots e_{j_1k_1}P_+\cdots)(\cdots)^\ast B^\ast_{j_1j_2}E_{-j_2k_1}(\cdots)^\ast,\nonumber
\end{eqnarray}
\begin{eqnarray}
&&\hspace*{-10mm}\sum_{j_1,j_2,k_1, k_2}\trace(\cdots B_{j_1j_2}B_{k_1k_2}\cdots)(\cdots B_{j_1j_2}B_{k_1k_2}\cdots)^\ast\label{gleichung_g}\\
&&=
\sum_{j_1,j_2,k_1}\bigl[\trace(\cdots e_{j_1k_1}K_1\cdots)(\cdots)^\ast B_{j_1j_2}B_{j_2k_1}(\cdots)^\ast\nonumber\\
&&\hspace*{15mm} +\trace(\cdots e_{j_1k_1}K_2 \cdots)(\cdots)^\ast B^\ast_{j_1j_2}B^\ast_{j_2k_1}(\cdots)^\ast\bigr]\nonumber.
\end{eqnarray}
From equations (\ref{gleichung_a}) to (\ref{gleichung_g}) we get the following conditions for the occuring summands
\begin{itemize}
 \item [(i')]
$E_{+j_1j_2}$ is followed only by $E_{-j_2k_1}$ or $B^\ast_{j_2k_1}$.
 \item [(ii')]
$E_{-j_1j_2}$ is followed only by $E_{+j_2k_1}$ or $B_{j_2k_1}$.
 \item [(iii')]
$B_{j_1j_2}$ is followed only by $E_{+j_2k_1}$ or $B_{j_2k_1}$.
 \item [(iv')]
$B^\ast_{j_1j_2}$ is followed only by $E_{-j_2k_1}$ or $B^\ast_{j_2k_1}$.
\end{itemize}
These conditions correspond to the former (i) to (iv).
Now let us sum over $j_1,\dots,j_{2r}$, $k_1\dots,k_{2r}$. 
\begin{eqnarray}
&& :=\!\!\!\!\!\!\!\sum_{j_1,\dots,j_{2r},k_1\dots,k_{2r}}\!\!\!\!\!\!\trace(X_{j_1j_2}^{(1)}X_{k_1k_2}^{(2)}\cdots X_{j_{2r-1}j_{2r}}^{(2r-1)}X_{k_{2r-1}k_{2r}}^{(2r)})(X_{j_1j_2}^{(1)}\cdots X_{k_{2r-1}k_{2r}}^{(2r)})^\ast\label{dicke_Summe}\\
&&=
(-1)^{L(w)}\!\!\!\!\!\!\!\!\!\!\!\!\sum_{j_1,\dots,j_{2r},k_1\dots,k_{r}}\!\!\!\!\!\!\!\!\!\trace(e_{j_1k_1}e_{j_3k_2}\cdots e_{j_{2r-1}k_r}\sigma(w))\tilde X_{j_1j_2}^{(1)}\tilde X_{j_2k_1}^{(2)}\!\!\!\!\cdots\tilde X_{j_{2r-1}j_{2r}}^{(2r-1)}\tilde X_{j_{2r}k_r}^{(2r)}\nonumber\\
&&=
(-1)^{L(w)}\!\!\!\!\sum_{j_1,\dots,j_{2r}}\!\!\trace(e_{j_1j_1}\sigma(w))\tilde X_{j_1j_2}^{(1)}\tilde X_{j_2j_3}^{(2)}\cdots \tilde X_{j_{2r}j_1}^{(2r)}\nonumber,
\end{eqnarray}
where $\tilde E_\pm=E_\mp$, $\tilde B=B^\ast$, $\tilde B^\ast=B$. Depending only on the word $w=\tilde X^{(1)}\cdots\tilde X^{(2r)}$, $\tilde X^{(l)}\in\{E_+,E_-,B,B^\ast\}$, there is some sign $(-1)^{L(w)}$
and a matrix $\sigma(w)$  which is a product of $r$ matrices of the form $K_{1/2}, P_\pm$.
In this way we get the sum over one type of word $w$ satisfying conditions (i)--(iv). All other words do not occur in $D_{2r}$.

We evaluate (\ref{dicke_Summe}) further. First we assume that $w$ is a word in $B,B^\ast$ only. Then $(-1)^{L(w)}=1$ and $\sigma(w)$ is a product in $K_1$ and $K_2$. As $K_1K_2=0=K_2K_1$, we get
\begin{eqnarray*}
 (\ref{dicke_Summe}) &=& \sum_{j_1,\dots,j_{2r}}\bigl(\trace(e_{j_1j_1}K_1)B_{j_1j_2}\dots B_{j_{2r}j_{1}}+
\trace(e_{j_1j_1}K_2)B^\ast_{j_1j_2}\dots B^\ast_{j_{2r}j_{1}}\bigr)\\
&=& \sum_{j_1,\dots,j_{2r}}\bigl(B_{j_1j_2}\dots B_{j_{2r-1}j_{2r}}+B^\ast_{j_1j_2}\dots B^\ast_{j_{2r}j_{1}}\bigr)\\
&=&\trace(B^{2r})+\trace((B^\ast)^{2r}).
\end{eqnarray*}
Now we allow $E_\pm$ to occur in $w$. Then $\sigma(w)$ is a product of $K_{1/2}$ and $P_\pm$. Notice that $P_+P_-=K_2$, $P_-P_+=K_1$, $P_+K_2=0=K_1P_+$, $P_+K_1=P_+=K_2P_+$, $P_-K_1=0=K_2P_-$ and $P_-K_2=P_-=K_1P_-$.
Thus, $\sigma(w)$ does not vanish if and only if $w$ satisfies (i) to (iv). Additionally, $\trace(\sigma(w))\not=0$ if and only if $P_+$ occurs exactly as often as $P_-$, i.e. if and only if
\begin{itemize}
 \item [(v)] $E_+$ occurs with the same multiplicity as $E_-$.
\end{itemize}
In this case $\trace(e_{j_1j_1}\sigma(w))=1$ and we have
\begin{equation*}
 (\ref{dicke_Summe}) =(-1)^{L(w)}\trace(w).
\end{equation*}
Thus, part (a) of the theorem  is proved apart from the sign $(-1)^{L(w)}$. To compute this sign, we must count the signs $(-1)$ given by equations (\ref{gleichung_c}) and (\ref{gleichung_e}) in the right way. That is, we find $L(w)$ to be the number of times $E_-B$ and $BE_+$ occur isolatedly in $w$ cyclicly.

For part (b) we notice that as long as $r<m$, the element $D_{2(r+1)}$ is not a polynomial in $D_2,\dots,D_{2r}$. For example, we never get $\trace((E_+E_-)^{r+1})$ as a combination of $\trace(E_+E_-),\dots,\trace((E_+E_-)^{r})$.
On the other hand it is well known and due to the Harish-Chandra isomorphism  that $\mathfrak Z(\mathfrak g^\CC)$ is generated by $m$ elements of length $2,\dots,2m$. So $D_2,\dots,D_{2m}$ must do.
\end{proof}



\section{Applications}

Let us assume we have an admissible representation $\Pi$ of $G=\Sp_m(\RR)$ and let us look at its isotypical component $\Pi_\rho$ for some irreducible representation $\rho$ of the maximal compact subgroup $K=K_m$. As $K$ is isomorphic to the unitary group $U_m$ by
\begin{equation*}
 J :K\to U_m, \quad\begin{pmatrix} A&-S\\S&A\end{pmatrix}\mapsto A+iS,
\end{equation*}
$\rho$ is characterized by its highest weight $(\lambda_1,\dots,\lambda_m)$. Let $v_h\not= 0$  be a highest weight vector of $\rho$. For $j\geq k$, the action of $B_{jk}$ on $v_h$ is determined by
\begin{eqnarray*}
 \rho(B_{jk})v_h &=&
\frac{d}{dt}\rho(\exp(tJ(B_{jk})))v_h\mid_{t=0}\\
&=&
\left\{
\begin{array}{ll}
 \frac{d}{dt}\rho(diag(1,\dots,e^{-t},1,\dots 1))v_h\mid_{t=0}=-\lambda_jv_h,&\textrm{ for }j=k,\\
\frac{d}{dt}\rho(\mathbf 1-te_{kj})v_h\mid_{t=0}=0,&\textrm{ for } j>k,
\end{array}\right.
\end{eqnarray*}
as $\exp(tJ(B_{jk}))$ is an upper triangular matrix if $j\leq k$.
Similarly we get for a lowest weight vector $v_l$,
\begin{eqnarray*}
 \rho(B_{kj})v_l &=&
\left\{
\begin{array}{ll} -\lambda_jv_l,&\textrm{ for }j=k,\\
0,&\textrm{ for } j>k.
\end{array}\right.
\end{eqnarray*}
Next we notice that for all words $w$ occuring in Theorem~\ref{Satz_Casimir_allgemein}, $\trace(w)$ is $\mathfrak k^\CC$-invariant (as we get telescopic sums for the commutators). Thus by Schur's lemma, the Casimirs' action on $\Pi_\rho$ is deduced by the actions of their single summands $\trace(w)$ on each $K$-irreducible component.On these components the summands are constant given by evaluating on the highest weight vector, for example.

Furthermore, $\trace(B^{2r})$, $\trace((B^\ast)^{2r})$ belong to $\mathfrak Z(\mathfrak k^\CC)$, so they act by constants on $\Pi_\rho$  deducible by $\rho(B_{jk})v_h$, $j\geq k$. For example, if we rearrange
\begin{equation*}
 \trace(BB)=\trace(B^\ast B^\ast)=\sum_j B_{jj}^2+\sum_{k<j}\bigl(2B_{kj}B_{jk}+B_{jj}-B_{kk}\bigr),
\end{equation*}
then we get
\begin{equation*}
 \rho(\trace(BB))v_h=\sum_j\bigl(\lambda_j^2+(m+1-2j)\lambda_j\bigr)v_h.
\end{equation*}
For the general case notice that in any summand $B_{j_1j_2}\dots B_{j_{2r}j_1}$ of $\trace(B^{2r})$ there is some $B_{j_nj_{n+1}}$ where $j_n>j_{n+1}$, if not all $j_n$ are equal. So by rearranging, we can determine the action of this summand by the action of terms of lower length.

For words $w$ in which both $E_\pm$ and $B,B^\ast$ occur, the evaluation of $\trace(w)$ is not that simple. But rearranging $\trace(w)$ (thereby producing terms of lower length satisfying again conditions (i)--(v) above) such that all terms $B, B^\ast$ are collected on the right, they can be evaluated first. For example, for the first two Casimirs (see Ex.~\ref{Prop_Casimir_2}) we get
\begin{cor}\label{cor_casimir_operatoren}
Let $C_1:=\frac{1}{2}D_2$ and $C_2:=\frac{1}{2}D_4$. Then
\begin{equation*}
 C_1=\frac{1}{2}(\trace(E_+E_-)+\trace(E_-E_+))+\trace(BB),
\end{equation*}
\begin{eqnarray*}
 C_2&=& \frac{1}{2}\bigl(\trace(E_+E_-E_+E_-)+\trace(E_-E_+E_-E_+)+\trace(B^4)+\trace((B^\ast)^4)\bigr)\\
&&+2\bigl(\trace(E_+E_-BB)+\trace(E_-E_+B^\ast B^\ast)\bigr)
-\sum_{i,j,k,l}\{(E_+)_{kl},(E_-)_{ij}\}B_{jk}B_{il}\\
&&+\frac{(m+1)^2}{2}(\trace(E_+E_-)+\trace(E_-E_+)),
\end{eqnarray*}
where 
\begin{equation*}
 \frac{1}{2}(\trace(E_+E_-)+\trace(E_-E_+))=\trace(E_+E_-)-(m+1)\trace(B),
\end{equation*}
\begin{eqnarray*}
&&\hspace*{-10mm}\frac{1}{2}(\trace(E_+E_-E_+E_-)+\trace(E_-E_+E_-E_+))=\trace(E_+E_-E_+E_-)\\
&&-\frac{1}{2}(\trace(E_+E_-)+\trace(E_-E_+))\trace(B)-\frac{m+2}{2}\bigl(\trace(E_+E_-B)+\trace(E_-E_+B^\ast)\bigr).
\end{eqnarray*}
\end{cor}

In the case $\rho=(\lambda,\dots,\lambda)$, we now have transparent formulae at hand. Here $\rho$ has dimension one, so highest and lowest weight vectors coincide and we have $\rho(B_{jk})=-\lambda\delta_{jk}$. So the only terms of $B,B^\ast$ left are $B_{jj}$ which produce a common constant. For example,
\begin{equation*}
 \Pi_\rho(C_1)= \Pi_\rho(\trace(E_+E_-))+\lambda m(m+1+\lambda)
\end{equation*}
and
\begin{eqnarray*}
 \Pi_\rho(C_2)&=&\Pi_\rho(\trace(E_+E_-E_+E_-))+m\lambda^4\\
&&+((m+1)^2+2\lambda (m+1)+2\lambda^2)\bigl(\Pi_\rho(\trace(E_+E_-)+\lambda m(m+1)\bigr).
\end{eqnarray*}
Similarly it is evident that
\begin{equation*}
 \Pi_\rho(D_{2r})=2\Pi_\rho(\trace((E_+E_-)^r))+\Pi_\rho(P_{2r}),
\end{equation*}
where $P_{2r}$ is a polynomial in $\trace(E_+E_-),\dots,\trace((E_+E_-)^{r-1})$. So we get
\begin{cor}\label{cor_casimir_allg_auf_rho}
 On $\Pi_\rho$, $\rho=(\lambda,\dots,\lambda)$, the Casimir operators are exactly the polynomials in $\trace(E_+E_-),\dots,\trace((E_+E_-)^{m})$.
\end{cor}
As an application, we consider modular forms on the Siegel halfplane $\H_m$.
For an irreducible  representation $(V,\rho)$ of $\GL_m(\CC)$ (equivalently of $K_m$),  let $f:\H_m\to V$ be a $C^\infty$-function  of moderate growth satisfying
\begin{equation*}
 f(g.z)=\rho(cz+d)f(z)
\end{equation*}
for all $g=\begin{pmatrix}a&b\\c&d\end{pmatrix}\in\Sp_m(\ZZ)$, $g.z=(az+b)(cz+d)^{-1}$. That is, $f$ is a non-holomorphic modular form for $\rho$. Then $f(g)=\rho^\ast(ci+d)f(g.i\mathbf 1_m)$ defines an automorphic form on $G$.
If more precisely $f$ is a modular form of weight $\kappa$, then $\rho^\ast=(-\kappa,\dots,-\kappa)$. By Corollary~\ref{cor_casimir_allg_auf_rho}, the action of $\mathfrak Z(\mathfrak g^\CC)$ on such modular forms is given by evaluating $\trace(E_+E_-),\dots,\trace((E_+E_-)^{m})$, which are differential operators on $\H_m$ (\cite{weissauersLN}, Ch.~3,~4). Especially, if $f$ is holomorphic, then $\trace((E_+E_-)^{r})f=0$.


\vspace*{.5cm}
\end{document}